\documentclass[english,a4paper,12pt,draft]{article}
\usepackage{babel,amsmath,amssymb,amsbsy,amsfonts,latexsym, amsthm}
\newtheorem{teor}{\bf Theorem}[section]

\newtheorem{cor}[teor]{\bf Corollary}
\newtheorem{rmk}{\bf Remark}

\newtheorem{lem}{Lemma}[section]
\newtheorem{pro}{Proposition}[section]

\theoremstyle{definition}

\newcommand{\Z}{\mathbb{Z}}
\newcommand{\R}{\mathbb{R}}
\newcommand{\N}{\mathbb{N}}

\begin{document}

\parindent 0pc
\parskip 6pt
\overfullrule=0pt

\title{\textbf{Moser-Trudinger inequality on conformal discs}}
\author{G. Mancini\thanks{Dipartimento di Matematica,
Universit\`a degli Studi "Roma Tre", Largo S. Leonardo Murialdo, 1 - 00146
Roma, Italy. E-mail  {\sf mancini@mat.uniroma3.it}.} 
and K. Sandeep\thanks{TIFR Centre for Applicable Mathematics, Sharadanagar,Chikkabommasandra, Bangalore 560 065. E-mail  
{\sf sandeep@math.tifrbng.res.in}}}

\date{}
\maketitle
\begin{abstract}
\begin{center}
\noindent We prove that a sharp Moser-Trudinger inequality holds true on a  conformal disc if and only if  the metric is bounded from above by the Poincar\'e metric. We also derive  necessary and sufficient conditions for the validity of a sharp Moser Trudinger inequality on a simply connected domain in $\R^2$
\end{center}
\end{abstract}

\section{Introduction} \label{intro}
 
In 1971 Moser, sharpening an  inequality due to Trudinger, proved that
\begin{equation} \label{MT}
 \sup\limits_{u\in H^1_0(\Omega), \quad  \int\limits_{\Omega}|\nabla u|^2 \leq 1}   \int\limits_\Omega \left(e^{4\pi u^2} - 1 \right) dx \: < \: +\infty
\end{equation} 
for every bounded open domain $\Omega \subset  \R^2$ (in \cite{M}). This inequality is sharp, in the sense that the 'critical' constant $4\pi$ cannot be improved. 
Referred as 'Moser-Trudinger inequality', (\ref{MT})  also implies the  estimate 
\begin{equation} \label{MTu}
\ln \frac{1}{|\Omega|}\int\limits_\Omega e^u \leq  C + \frac{1}{16 \pi} \int\limits_{\Omega}|\nabla u|^2 \qquad \forall u\in H^1_0(\Omega)
\end{equation}
for some universal constant $C>0$ and any bounded domain $\Omega$ and, again, the constant  $\frac{1}{16 \pi}$ is sharp (see \cite{KL}). 
\\ Inequality (\ref{MT}) has been extended to any 2-d  compact Riemannian manifold with or without boundary (see \cite{B}, \cite{F}) or even to some subriemannian manifolds (see \cite{BFM} and references therein). However, little is known in case $\Omega$ is a non compact 2-d Riemannian manifold, even in the simplest cases $\Omega \subset \R^2$ with $|\Omega|=\infty$ (see \cite{AT}, \cite{RU}) or $\Omega=\mathcal H^2$, the 2-d hyperbolic space.\\
We address here the case of conformal discs, i.e. $\Omega=D$, the unit open disc  in $\R^2$, endowed with a conformal metric $g= \rho g_e$, where $g_e$ denotes the euclidean metric and $\rho\in C^2(D), \: \rho > 0$. 
Denoted by $dV_g= \rho dx$ the  volume form, by conformal invariance of the Dirichlet integral (\ref{MT}) takes the form
\begin{equation} \label{MTc}
\sup\limits_{u\in C^\infty_0(D), \quad \int\limits_{D}|\nabla u|^2 \leq 1}   \int\limits_D \left(e^{4\pi u^2} - 1 \right) dV_g < \infty
\end{equation}
A relevant case is  the hyperbolic metric $g_h := (\frac{2}{1-|x|^2})^2 g_e$. We will show that (\ref{MTc}) holds true in this case. Actually, we have the following

{}{}{}{}{}{}{}{}{}{}{}{}

\begin{teor} \label{conf}
Given a conformal metric $g$ on the disc, (\ref{MTc}) holds true if and only if $g \leq c g_h$ for some positive constant $c$. 
\end{teor}
After a personal communication,in \cite{ATi} the inequality \eqref{MTc} with $g=g_h $ found an application in the study of blow up analysis and eventually a different proof of \eqref{MTc} when $g=g_h $.\\\\
As for (\ref{MT}) in case $|\Omega|=+\infty$, the supremum therein will be in general infinite. To have it finite, an obvious  necessary condition  is that
$$ \lambda_1(\Omega) := \inf \left\{\int\limits_\Omega |\nabla u |^2 dx: \quad u\in C^\infty_0(\Omega), \int\limits_\Omega |u|^2 = 1 \right\} > 0$$
As a partial converse, it was shown by D.M. Cao \cite{C}   that $\lambda_1(\Omega)>0$ implies subcritical exponential integrability, i.e. for every $\alpha < 4 \pi$ it results
\begin{equation} \label{MTsc}
 \sup\limits_{u\in C^\infty_0(\Omega), \quad  \int\limits_{\Omega}|\nabla u|^2 \leq 1}   \int\limits_\Omega \left(e^{\alpha u^2} - 1 \right) dx < \infty
\end{equation}
(see also \cite{O} and \cite{AT}, for a scale invariant version of Trudinger 
inequality which implies (\ref{MTsc}) ).
However, \emph{no information is provided for  the  critical case $\alpha = 4\pi$}. We will show that  $\lambda_1(\Omega)>0$ is, on simply connected domains, also sufficient for   (\ref{MT}) to hold true. To state our result, let
$$\omega (\Omega):= \sup \{r>0 :\quad \exists  D_r(x) \subset \Omega \}$$
\begin{teor} \label{width}
Let $\Omega$ be a simply connected domain in $\R^2$. Then 
\begin{center}
 (\ref{MT}) holds true \quad$\Leftrightarrow$ \quad
$\lambda_1(\Omega) > 0$\quad$\Leftrightarrow$ \quad
$\omega (\Omega) < +\infty$
\end{center}
\end{teor}
\begin{rmk}
The topological assumption on $\Omega$ cannot be  dropped: in Appendix we exhibit  domains  $\Omega$ with 
$\omega (\Omega)<+\infty$ and $\lambda_1(\Omega)=0$, for which, henceforth,  (\ref{MT}) fails.
However, we suspect that $\lambda_1(\Omega)>0$ is sufficient to insure (\ref{MT}).
\end{rmk}
\setcounter{equation}{0}

\section{Proof of the main results and  asymptotics for $L^p$ Sobolev inequalities} \label{Hs1}

\textbf{Proof of Theorem \ref{conf}}.
Let us write $g := \rho g_e = \zeta g_h$ whith $ \zeta := \rho \frac{(1-|x|^2)^2}{4}$. \\
We first prove that if there is $x_n\in D$ such that $\zeta(x_n)\rightarrow_n +\infty$, then there are $u_n \in H^1_0(D)$ with $\int\limits_{D}|\nabla u_n|^2 \leq 1$ such that $\int\limits_D \left(e^{4\pi u^2_n} - 1 \right) \zeta dV_h \rightarrow_n +\infty $. 
To this extent, let $\varphi_n$ be a conformal diffeomorphism of the disc such that $\varphi_n(0)=x_n$. Then
$$\exists \epsilon_n > 0 \quad \textit{such that} \qquad  |x|\leq \epsilon_n \quad \Rightarrow \quad \zeta(\varphi_n(x)) \geq \frac 12  \zeta(\varphi_n(0))\rightarrow_n +\infty$$
Let now $v_n(x)=v_n(|x|)$ be the Moser function defined as
$$v_n(r) \: = \:  \sqrt{\frac{1}{2\pi}} \:\left[ \left(\log\frac{1}{\epsilon_n}\right)^{\frac 12} \chi_{[0, \epsilon_n)} \:  
+  \:\left( \: \log\frac{1}{\epsilon_n}\right)^{-\frac{1}{2}}\:  \log \frac 1r \: \: \:  \chi_{[\epsilon_n, 1]}\right] $$ 
Notice that $\int\limits_D |\nabla v_n|^2 = 1$. Let $u_n := v_n \circ \varphi_n^{-1}$. Then, by conformal invariance  and because $\varphi_n$ are hyperbolic isometries,
$\int\limits_D |\nabla u_n|^2 = 1$ and $$\int\limits_D \left(e^{4\pi u^2_n} - 1 \right) \zeta dV_h = \int\limits_D \left(e^{4\pi v^2_n} - 1 \right) \zeta \circ \varphi_n dV_h $$
$$\geq
\frac{\zeta(x_n)}{2} \: \left(\frac{1}{\epsilon_n^2} - 1\right) \int\limits_{|x|\leq \epsilon_n} dV_h \: = \: \frac{\zeta(x_n)}{2} \: \frac{1-\epsilon_n^2}{\epsilon_n^2} \: \frac{4\pi \epsilon_n^2}{1-\epsilon_n^2} \rightarrow_n +\infty$$
Hence  a bound for $g$ in terms of $g_h$ is necessary for (\ref{conf}) to hold true.\\
We now prove that boundedness of $\zeta$ is also sufficient
for (\ref{conf}) to hold true. Under this assumption, (\ref{conf}) reduces to
\begin{equation} \label{MTT}
\sup\limits_{u\in C^\infty_0(D), \quad  \int\limits_{D}|\nabla u|^2 \leq 1}   \int\limits_D \left(e^{4\pi u^2} - 1 \right)   dV_h < \infty
\end{equation}
Let $u^\ast $  be  the symmetric decreasing hyperbolic rearrangement of $u$, i.e.
$$\mu_h(\{u^\ast > t \}) = \mu_h(\{u > t \}) $$  
By the properties of the rearrangement (see \cite{Bae}), it is enough to prove  (\ref{MTT}) for $u$ radially symmetric. For  $u$   radial, inequality (\ref{MTT}) rewrites, in hyperbolic polar coordinates $|x|= \tanh \frac t2$, as
\begin{equation} \label{MTTT}
\sup\limits_{ 2\pi \int\limits_0^\infty |u'|^2 \sinh t \: dt \: \leq 1}   \int\limits_0^\infty \left(e^{4\pi u^2} - 1 \right)  \: \sinh t \: dt < \infty
\end{equation}
To prove (\ref{MTTT}), observe first that from $\int\limits_{D}|\nabla u|^2 \leq 1$ it follows , for $t<\tau$,\\
$|u(\tau) - u(t)| = \left| \int\limits_t^{\tau} u'(s)  \: ds \right| \leq \left(\int_t^\infty |u'|^2 \sinh s \: ds \right)^{\frac 12} \times \left( \int\limits_t^\infty \frac { ds}{\sinh s}  \right)^{\frac 12} \: \leq \: \left( \frac{1}{2\pi \: \sinh t }\right)^{\frac 12}$
and since  $\int\limits_D u^2 dV_h = 2\pi \int\limits_0^\infty u^2 \sinh t \: dt < +\infty$
implies $\liminf\limits_{\tau\rightarrow +\infty} u(\tau) = 0$, we get
\begin{equation}\label{decayu}
\int\limits_{D}|\nabla u|^2 \leq 1 \quad \Rightarrow \quad |u(t)| \leq  \left(\frac{1}{2\pi \: \sinh t}\right) ^{\frac 12} \qquad \forall t
\end{equation}
Now, let  $ \: 2\pi \sinh T > 1 \: $  so that $\: \int\limits_{D}|\nabla u|^2 \leq 1 \:$ implies $ \: u(T)<1$,\quad and set\\
$v:= u - u(T), \: w:= \sqrt{1+u^2(T)} \:  v \: $ so that  $ \: w(T)=0 \: $ and  $ \: 2\pi \int\limits_0^T |w'|^2 t \: dt$\\
$ \leq \: 
 [1 + u^2(T)]  \: 2 \pi \: \int\limits_0^T |u'|^2 \sinh t  \: dt  \leq
\left[1 + \left(\int\limits_T^\infty |u'|^2  \sinh t \: dt \right) \: \left(\int\limits_T^\infty \frac{dt}{\sinh t} \right)\right]\times$\\
$  \left[1 -  \int\limits_T^\infty 2\pi|u'|^2 \sinh t \: dt \right] \leq
\left[1 + \left(\int\limits_T^\infty |u'|^2  \frac{\sinh t}{\sinh T} \: dt \right) \:   \right] \times \left[1 -  \int\limits_T^\infty 2\pi|u'|^2 \sinh t \: dt \right] $\\
$\leq \left[1 + 2\pi \int\limits_T^\infty |u'|^2  \sinh t \: dt  \right] \times \left[1 - 2\pi \int\limits_T^\infty |u'|^2 \sinh t \: dt \right]\leq 1$.\\
Now, an application of  (\ref{MT}) gives
\qquad $2\pi \int\limits_0^T e^{4\pi w^2} \: t \: dt \leq c T^2 $ \qquad 
and, since\\  $u^2 = v^2 + 2 v u(T) + u^2(T) \leq v^2 + v^2 u(T)^2 +1+ u^2(T) \leq w^2 + 2$ \qquad implies\\
$2\pi \int\limits_0^T e^{4\pi u^2}    \sinh t \: dt \leq 2e^{8\pi} \: \pi  \:  \frac{\sinh T}{T} \: \int\limits_0^T e^{4\pi w^2} t \:dt
\leq c(T) T^2 $ , \quad we get
\begin{equation} \label{bo}
2\pi \int\limits_0^T [e^{4\pi u^2}-1]    \: \sinh t \: dt 
\: \leq \:\: c(T)
\end{equation}
for some constant $c(T)$ which does not depend on $u$. \\ 
Now, using (\ref{decayu}) and Hardy inequality \quad $\int\limits_D |\nabla u|^2 \geq \frac 14 \int\limits_D |u|^2 dV_h$, \quad we get\\
$\int\limits_T^\infty [e^{4\pi u^2}-1]    \: \sinh t \: dt \leq      \left[2 \int\limits_D |u|^2 \: dV_h +  
\sum_2^\infty \frac{(4\pi)^p}{p!} \int\limits_T^\infty |u|^{2p} \sinh t \: dt\right] \leq$\\ 
$\left[8 \int\limits_D |\nabla u|^2 \: dV_h +  \sum_2^\infty \frac{(4\pi)^p}{p!} \int\limits_T^\infty [\frac{1}{2\pi \sinh t}]^p \sinh t \: dt \right]
\left[8 +  \sum_2^\infty \frac{2^p}{p!} \int\limits_T^\infty \frac {dt}{(\sinh t)^{p-1}}\right]$.\\
 From
$\int\limits_T^\infty \frac {dt}{(\sinh t)^{p-1}}=\int\limits_T^\infty \left[\frac{2}{e^t-e^{-t}}\right]^{p-1} \: dt =
2^{p-1} \int\limits_T^\infty \frac{e^{-(p-1)t}}{(1-e^{-2t})^{p-1}} \: dt $ 
$\leq\left[\frac{2}{1-e^{-2T}}\right]^{p-1} \: \frac{e^{-(p-1)T}}{p-1}$\\
$ = \frac{1}{p-1} \: \left[\frac{1}{\sinh T}\right]^{p-1} \leq \left[\frac{1}{\sinh T}\right]^{p-1}$ if $p\geq 2$ \qquad and the  above inequality we get
\begin{equation} \label{boo}
2\pi \: \int\limits_T^\infty [e^{4\pi u^2}-1] \:    \: \sinh t \: dt \: \leq \:  2 \pi \:   \: 
\left[ 8 + \sinh T \: \:  e^{\frac{2}{\sinh T}} \right]\: = \: c(T)
\end{equation}
Inequalities (\ref{bo}) and (\ref{boo}) give (\ref{MTTT}) and hence (\ref{conf}).\\
 \textbf{Proof of Theorem \ref{width}  }  (\ref{MT})   implies $4\pi \int\limits_\Omega (\frac{u}{\|\nabla u\|})^2 \: dx \leq \int\limits_\Omega \left(e^{4\pi (\frac{u}{\|\nabla u\|})^2} - 1 \right) dx$\\$ \leq c(\Omega)$ and hence 
 $\lambda_1(\Omega) \geq  \frac{4\pi}{c(\Omega)}$. In turn, this clearly implies
 $\omega (\Omega) < +\infty$.\\
To complete the proof, it remains to show that if $\Omega$ is simply connected then
$\omega (\Omega) < + \infty $ implies (\ref{MT}).
Let $\varphi : D \rightarrow \Omega$ be a conformal diffeomorphism, so that 
(\ref{MT}) rewrites as (\ref{MTc})
where $g:= \varphi^\ast g_e = |\det J_\varphi| g_e$. Let us show that
\begin{equation} \label{K}
\omega (\Omega) < R \quad \Rightarrow \quad  |det J_\varphi (x)| \leq  \frac{16 R^2}{(1-|x|^2)^2}
\end{equation} 
so that Theorem \ref{conf} applies to give the conclusion. Now, (\ref{K}) follows from \\
Koebe's covering Theorem (see \cite{G}): \emph{if  $\psi : D \rightarrow \Omega$  is a conformal diffeomorphism and $z \notin \psi(D)$ for some $z\in D_r(\psi(0))$ , then $|\psi'(0)| \leq 4r$}. In fact, given $w\in D$, let $ \varphi_w (z):= \varphi (w + (1-|w|) z), \quad z\in D$. 
By assumption,  $\varphi_w (D)= \varphi (D_{1-|w|}(w))$ cannot cover the disc $D_R(\varphi(w))$, and hence  $|det J_\varphi (w)|^{\frac 12} |(1-|w| ) |= |\varphi'(w) \: (1-|w|) | = | \varphi_w'(0)| \leq 4 R$.

We end this Section  deriving, from Moser-Trudinger inequalities, an asymptotic formula for best constants in $L^p$ Sobolev inequalities on 2-d Riemannian manifolds $(M,g)$ (see \cite{RW}, \cite{AG}, for smooth bounded domains in $\R^2$).
For notational convenience, we say that  $(M,g)$ is an MT-manifold if 
\begin{equation}\label{MTM}
\sup\limits_{u\in C^\infty_0(M), \quad \int\limits_{M}|\nabla_g u|^2 dV_g \leq 1}   \int\limits_M \left(e^{4\pi u^2} - 1 \right)  dV_g < \infty
\end{equation} 
\begin{pro}\label{asy}
Let $(M,g)$ be  an MT-manifold. Then 
\begin{equation}\label{asym}
S_p \: = \: S_p(M,g)\: := \:  \inf_{u\in C_0^\infty (M), u\neq 0} \frac{\int\limits_M |\nabla_g u|^2  dV_g}{(\int\limits_M |u|^p  dV_g)^{\frac 2p}} \: = \: \frac{8\pi e + \circ(1)}{p} 
\end{equation}
\end{pro}
\begin{proof}
 Let us prove first
\begin{equation}\label{liminf}
\liminf_p  \: p S_p \: \geq \:  8\pi e 
\end{equation}
By assumption, there is $C > 0$ such that, for every $p \in \N$, it results
$$\int\limits_M |\nabla_g u|^2 \:  dV_g\leq 1\quad \Rightarrow \quad C \geq \int\limits_M (e^{4\pi u^2} - 1) \:  dV_g  \geq  \frac{(4\pi)^p}{p!} \int\limits_M  |u|^{2p} \:  dV_g $$
and hence 
$$  \left(\int\limits_M  |u|^{2p} \:  dV_g \right)^{ \frac{1}{2p}}\leq \frac{C^{\frac{1}{2p}} (p!)^{\frac{1}{2p}}}{\sqrt{4\pi}} 
(\int\limits_M |\nabla_g u|^2  dV_g)^{\frac 12} \quad \forall u\in C_0^\infty(D)$$
If $n\leq p \leq n+1$, let $\alpha = \frac{n(n+1-p)}{p}$ and get, by interpolation,\\

 $\|u\|_{2p} \leq \frac{1}{\sqrt{4\pi}}C^{\frac{\alpha}{2n}}(n)!^{\frac{\alpha}{2n}}\times
 C^{\frac{1-\alpha}{2(n+1)}}(n+1)!^{\frac{1-\alpha}{2(n+1)}} \|\nabla_g u\| \qquad
 $ \qquad
 and hence\\
  
 $S_{2p} \geq \frac{4\pi}{C^{\frac 1p} (n!)^{\frac 1p} (n+1)^{\frac{1-\alpha}{n+1}}}
 \geq \frac{4\pi }{C^{\frac 1p} (n!)^{\frac 1n} (n+1)^{1-\frac{n}{p}}}$. \quad
 By Stirling's formula we obtain\\
 
$2p S_{2p} \geq \frac{8 p \pi e}{C^{\frac 1p} n [(1+\circ(1)) \sqrt{2\pi n}]^{\frac 1n} (n+1)^{1-\frac{n}{p}}}\geq
\frac{8  \pi e}{1+\circ(1)}$ \qquad and hence (\ref{liminf}).\\

To prove the reverse inequality, we use again the Moser function. For fixed $R>0$ and $0<l<R$, define $M_{l}(x)=M_l(|x|)$ on $\R^2$ as follows:
$$ M_{l}(r) = \sqrt{\log (\frac Rl)} \: \left[ \chi_{[0,l)}+\frac{\log (\frac {R}{r})}{\log (\frac Rl)}\chi_{[l, R)} \right], \qquad r\geq 0$$
Let $q\in M$ and choose $R>0$  strictly less than the injectivity radius of $M$ at $q$ and define $u_l(z) := M_{l}(Exp^{-1}_q (z))$ where $Exp_q$ is the exponential map at $q$. Note that $u_l$ is well defined and in $H^1(M).$
Now calculating in normal coordinates we get $$\int\limits_M |\nabla_g u_l|^2\;dV_g = \int\limits_{B(0,R)}g^{i,j}(x)(M_l)_{x_i}(M_l)_{x_j} \sqrt{g(x)} dx $$
Since the metric is smooth and $g_{i,j}(0) = \delta_{i,j}$ we get  $g^{i,j} = \delta_{i,j} + O(|x|)$ and $\sqrt{g(x)} = 1 + O(|x|)$. Using this we get
$$\int\limits_M |\nabla_g u_l|^2\;dV_g = 2\pi + O(1) (\log \frac Rl)^{-1}$$
Similarly
$$\int\limits_M |u_l|^{p}\;dV_g \ge \int\limits_{B(0,l)}|M_l(x)|^{p} \sqrt{g(x)} dx
= C(\log \frac Rl)^{\frac{p}{2}} l^2 $$ for some $C>0.$
Taking $\log\frac Rl = \frac{p}{4} $ and sending $p$ to infinity, we get
$$\limsup\limits_{p\rightarrow \infty} \:  p S_p \: \le \: \lim\limits_{p\rightarrow \infty}\frac{ \int\limits_M |\nabla_g u_l|^2\;dV_g }{\left(\int\limits_M |u_l|^{p}\;dV_g \right)^{\frac{2}{p}} } \: \le \: 8\pi e$$
\end{proof} 
\begin{cor}\label{asym}
If $g \leq  c g_h$ then $S_p(D, g) = \frac{8\pi e + \circ(1)}{p}$ 
\end{cor} 
\begin{rmk}\label{bad}
Let $p\in [1, 2)$ and $u_p= (1-|x|^2)^{\frac 1p} $. Then  \quad  $u_p\in H^1_0(D)$ \quad  
and
\quad $ \int\limits_D |u_p|^p dV_h = 
\int\limits_D |e^u - 1| dV_h = + \infty$. \quad In particular, $S_p=0$ for $p\in [1, 2)$.
\end{rmk}
Let us now derive from Proposition \ref{asy} an inequality analogous to  (\ref{MTu}).
\begin{cor}\label{MTO}
Let $(M,g)$ be an MT-manifold.  Then, chosen $  \delta\in (0,1)$, there is a constant $C(\delta)>0$ such that, for every  $u\in H^1_0(M)$,   it results  
\begin{equation} \label{MTOD}
\ln \int\limits_{M} [e^{u} -1]^2  dV_g \leq \ln \int\limits_{M} [e^{2u} -2u - 1]  dV_g \leq C(\delta) + \frac{1}{4 \delta \pi} \int\limits_{M}|\nabla_g u|^2 d V_g 
\end{equation} 
\end{cor}
\begin{proof}
After fixing $\delta\in (0,1)$, we get, by Taylor expansion 
$$\int\limits_{M} [e^{u} -u - 1]  dV_g 
= \sum\limits_{p=2}^\infty \frac{1}{p!}\int\limits_M u^p dV_g  \leq \sum\limits_{p=2}^\infty \frac{1}{\sqrt{p!}} \left[\frac{\|\nabla_g u\|^2}{8\pi \delta}\right]^{\frac p2} \left(\frac{8\pi \delta}{S_p}\right)^{\frac p2}\frac{1}{\sqrt{p!}} $$ 
$$\leq \left[\sum\limits_{p=2}^\infty \frac{1}{p!} \left(\frac{\|\nabla_g u \|^2}{8\pi \delta}\right)^p \right]^{\frac 12} \: \times \: \left[\sum\limits_{p=2}^\infty \frac{1}{p!}  \left(\frac{8\pi \delta}{S_p}\right)^p\right]^{\frac 12}$$
Since, by Stirling's formula and   (\ref{liminf}) 
\qquad $\limsup_p \frac{1}{(p!)^{\frac 1p}}  \frac{8\pi \delta}{S_p} \leq \delta < 1$ \qquad we\\
 
conclude, also using the inequality  $(e^t-1)^2 \leq
e^{2t} -2t - 1, \quad \forall t\in \R$, that
\begin{equation} \label{MTODa}
\int\limits_{M} [e^{u} -1]^2  dV_g \leq \int\limits_{M} [e^{2u} -2u - 1]  dV_g  \leq  c(\delta) \left(e^{\frac{\|\nabla_g u \|^2}{2\pi \delta}}-\frac{\|\nabla_g u \|^2}{2\pi \delta} - 1\right)^{\frac 12}
\end{equation}   
\end{proof}
\begin{rmk} \label{omega} We believe that  (\ref{MTOD}) holds with $\delta=1$ (and $\frac{1}{4  \pi}$ is optimal). Actually, as it is clear from the proof,  subcritical exponential integrability  (\ref{MTsc}) is enough to get (\ref{MTOD}).  In particular,  (\ref{MTOD}) holds true if $M=\Omega$, a smooth open subset of $\R^2$ with $\lambda_1(\Omega)>0$.   
\end{rmk}

\setcounter{equation}{0}
\section{Application to a geometric PDE } \label{curv}
Here we apply Moser-Trudinger inequality  to the following problem.\\ 
\textit{ Let $\Omega$ be a smooth open set in $\R^2$. Let $K \in C^\infty(\Omega)$.\\Is it  $K$ the Gauss curvature of a conformal metric $g=\rho g_e$ in $\Omega$?}\\
It is known that solving this problem amounts to solve the equation
\begin{equation} \label{curvature zero}
\Delta v  + K e^{2v} = 0 \qquad \textit{in} \quad \Omega
\end{equation}
In fact, if  $v\in C^2(\Omega)$  solves (\ref{curvature zero}) then $e^{2v}g_e$ is a  conformal metric  having $K$ as Gauss curvature.
Equation (\ref{curvature zero}) is not  solvable in general, e.g. if $\Omega=\R^2$, $K \leq 0$ and $K(x) \leq -|x|^{-2}$ near $\infty$ (a result due to Sattinger, see \cite{CL} or \cite{KW}). 
In \cite{KW} it is also noticed, as a Corollary of a  general result, that if $\Omega$ is bounded and $K\in L^p(\Omega)$ for some $p>2$, then (\ref{curvature zero}) is solvable. We prove 
\begin{teor} \label{KW} Let $(M,g)$ be  an MT-manifold. Let $K_i \in L^2(M)$. Then equation
  \begin{equation} \label{Mc} 
\Delta_g v  + K_1 + K_2 e^{2v} = 0
\end{equation} 
has a  solution in $H_0^1(M)+\R$
\end{teor}
\begin{rmk} In view of Remark \ref{omega}, Theorem \ref{KW} applies to any smooth open set $\Omega \subset \R^2$ for which $\lambda_1(\Omega)>0$.
\end{rmk} 
When $\Omega$ is the unite disc, sharp existence/nonexistence results for (\ref{curvature zero}) have been obtained by Kalka and Yang \cite{KY} in the case of \textsl{nonpositive} $K$.
The following result is a restatement of Theorem 3.1 in \cite{KY}:
\begin{teor}(Kalka and Yang) \label{KY} 
Let $K\in C(D)$, $K < 0$ in $D$. Assume
 $$\exists \alpha > 1, \: C>0 \quad \textit{such that} \quad K  \: \geq \:  - \:  \frac{C }{(1-|x|^2)^2 | \log(1-|x|^2)|^{\alpha}} $$
Then equation (\ref{curvature zero}) has a $C^2$ solution.
 If $$K \: \leq  \: - \: \frac{C }{(1-|x|^2)^2 | \log(1-|x|^2)|}\qquad \textit{for $|x|$ close to $1$} $$ then 
 (\ref{curvature zero}) has no $C^2$ solution in $D$.
\end{teor}
Existence   is proved by monotone iteration techniques. 
We present here a variational existence result \textit{without sign assumptions on $K$}.  
\begin{teor} \label{var} 
Let   $\int\limits_D K^2 (1-|x|^2)^2 dx < +\infty$.
Then equation (\ref{curvature zero})  has a  solution in $H_0^1(D)+\R$.
\end{teor}
\begin{rmk} This result is far from being sharp. For instance, if one takes $K_\alpha=-\frac \alpha 2 (\frac{2}{1-|x|^2})^{2-\alpha}, \quad \alpha \in \R$, (\ref{curvature zero}) has the solution $v_\alpha = \frac \alpha 2 \log\frac{2}{1-|x|^2}$, so that $K_\alpha$ is the curvature of $g_\alpha = (\frac{2}{1-|x|^2})^\alpha g_e$. So, negative $\alpha $ give examples of \textsl{positive}
curvatures $K_\alpha$ with arbitrary blow up. 
\end{rmk}
Proofs of Theorems \ref{KW} and \ref{var} rely on inequality (\ref{MTODa}). We state below  some  consequences of (\ref{MTODa}) that we need. 
\begin{lem}\label{cont}
Let $(M,g)$ be an $MT$ manifold. Let $K \in L^2(\mu_g)$.
 Then\\ $I_K(v):= \int\limits_M K [e^{v} -1] dV_g$\quad is uniformly continuous on bounded sets of $ H^1_0(M)$.
Furthermore, \\ $v_n\rightharpoonup v$ in $H^1_0(M)$ implies $I_K(v_n) \rightarrow I_K(v)$.
\end{lem}
\begin{proof} Let $\|\nabla_g u\|+\|\nabla_g v\| \leq R$. Writing $e^t-e^s= (e^{t-s}-1) (e^s-1)+ (e^{t-s}-1)$, we see, using  the inequality $(e^t-1)^2 \leq |e^{2t}-1| \: \:  \forall t$ and (\ref{MTODa}),  that
$$|I_K(u)-I_K(v)| \leq (\int\limits_M K^2 dV_g )^{\frac 12} \times (\int\limits_M |e^{u} -e^v|^2 dV_g )^{\frac 12}\leq$$
$$c(K) \left[(\int\limits_M |e^{2v} -1|^2 dV_g )^{\frac 14} \times (\int\limits_M |e^{2(u-v)} -1|^2 dV_g )^{\frac 14} + (\int\limits_M |e^{u-v} -1|^2 dV_g )^{\frac 12}\right]\leq$$
$$c(K, R, \delta)  \left(e^{\frac{2\|\nabla_g (u-v) \|^2}{\pi \delta}}-\frac{2\|\nabla (u-v) \|^2}{\pi \delta} - 1\right)^{\frac 18}
\leq C(K, R, \delta)\|\nabla_g (u-v)\|^{\frac 12}$$
Next, assume $v_n\rightharpoonup v$ in $H^1_0(M)$ and a.e. From  $\sup_n \int\limits_M |\nabla_g v_n|^2 < \infty$ and Lemma \ref{MTO} we get
$\sup_n \int\limits_M  (e^{v_n}-1)^2 dV_g < +\infty$ and hence Vitali's convergence theorem
applies to get
$\int\limits_A K (e^{v_n}-1) dV_g \rightarrow_n \int\limits_A K (e^{v}-1) dV_g$.
\end{proof}
We state without proof the following property
\begin{cor}\label{lsc}
Let $(M,g)$ be an $MT$ manifold. Let   $I(v):= \int\limits_M [e^{v} -1]^2 dV_g$, $J(v):= \int\limits_M [e^{v} -v - 1] dV_g$. Then $I, J \in Lip_{loc}(H^1_0(M))$. 
\end{cor}

\textbf{Proof of Theorem \ref{KW}} \quad Let $O:=\{v\in H^1_0(M): \: \; \int\limits_M K_2 (e^{2v}-1)dV_g>0 \}$.  By Lemma \ref{cont}, $O$ is open. Let 
$$E_K(v)=  \int\limits_M |\nabla_g v|^2 dV_g - 2 \int\limits_M K_1 \:v \: dV_g - \log  \int\limits_M K_2 (e^{2v}-1) dV_g \qquad v\in O$$
Since (\ref{MTM}) implies $\lambda_1(g):= S_2(M,g) >0$, we get from  Corollary \ref{MTO} and the assumption on $K_i$,  $$E_K (v) \: \geq \:  \int\limits_M \: |\nabla_g v|^2 \: dV_g \: -$$ 
$-\left[\frac{2}{\lambda_1(M)}(\int\limits_M K_1^2  \: dV_g )^{\frac 12} (\int\limits_M  |\nabla_g v|^2  \: dV_g )^{\frac 12}  + c(K_2) +
 \frac 12   \log   \int\limits_M  (e^{2v}-1)^2 dV_g  \right]\geq $\\ 
$\geq \left(\int\limits_M |\nabla v|^2 dV_g\right)^{\frac 12}
\left[( 1-\frac{1}{\pi \delta})  \: \left(\int\limits_M |\nabla v|^2 dV_g\right)^{\frac 12} - c(K_1, M)\right]  \: - c(K_2, \delta)$\\

for every $ v\in O$. Thus $E_K$ is bounded below and coercive on $O$. Hence, if $v_n \in O$ ,  $ E_K(v_n) \rightarrow \inf\limits_{ O} E_K$, we can assume $v_n$ converges weakly to some $v$. By Lemma \ref{cont} and boundedness of $E_K(v_n)$ we infer that $v\in O$ and $E_K(v) = \inf\limits_{ O} E_K$. Since  $O$ is open, we see that
$$ \int\limits_M \left[\nabla_g v \nabla_g \varphi -  K_1 \varphi\right] dV_g -
\frac{\int\limits_M K_2 e^{2v}\varphi dV_g}{\int\limits_M K_2 (e^{2v}-1) dV_g}=0 \quad \forall \varphi \in C_0^\infty(M)$$
and hence $v-\frac 12 \log \int\limits_M K_2 (e^{2v}-1) dV_g$ solves (\ref{curvature zero}).\\

\textbf{Proof of Theorem \ref{var}} \quad It goes like above, with the obvious modification
$$E_K (v) :=  \int\limits_D |\nabla v|^2 dx  - \log  \int\limits_D K (e^{2v}-1) dx \geq $$  $$\int\limits_D |\nabla v|^2 dx - \log   (\int\limits_D K^2 (1-|x|^2)^2 dx )^{\frac 12} \;\:  (\int\limits_D  \frac{(e^{2v}-1)^2}{(1-|x|^2)^2} dx )^{\frac 12} \geq $$
$$ (1-\frac{1}{2\pi \delta})\int\limits_D |\nabla v|^2 dx - c(K, \delta) \qquad \forall v\in O$$

\begin{rmk}
 In  \cite{WL} a similar result is proven, but under the stronger assumption  $|K(x)| \leq \frac{C}{(1-|x|)^\alpha}$ with $\alpha \in (0,1)$. 
\end{rmk}
The result in Theorem \ref{var}, when applied to negative $K$, is weaker than the one in Kalka-Yang. But, even more, the solutions we find don't address the main point in
\cite{KY}, which  is to find \textsl{complete metrics} of prescribed (nonpositive) Gaussian 
curvature on noncompact Riemannian surfaces:  a solutions of (\ref{curvature zero}) has to blow to $+\infty$ along $\partial D$ to give rise to a complete metric, and this is not  the case for the solutions obtained in Theorem \ref{var}. A first step in this direction is to build solutions of (\ref{curvature zero}) with prescribed boundary values.\\ Since without sign assumptions on $K$ one cannot expect $K$ to be the curvature of a complete  metric $g$ ( e.g., if $K \geq 0$ around $\partial D$, then $K$ cannot be the curvature of a complete conformal metric on the disc (see \cite{KY}))  we restrict our attention to  $K<0$. Assuming again $\int\limits_D K^2 (1-|x|^2)^2 dx < +\infty$, we see that the strictly convex functional
 $$J_K(v)= \frac 12 \int\limits_D |\nabla v|^2 dx - \int\limits_D K  (e^{2v}-1)dx
 \qquad v\in H^1_0(D)$$
 is well defined, uniformly continuous and weakly lower semicontinuous by Lemma \ref{cont}. Furthermore, by Hardy's inequality,
 $$J_K(v)= \frac 12 \int\limits_D |\nabla v|^2 dx - \int\limits_D K  v dx  - \frac 12 \int\limits_D K (e^{2v}-2v-1) dx \geq $$
 $$\frac 12 \int\limits_D |\nabla v|^2 dx -
\frac{1}{2}(\int\limits_D K^2 (1-|x|^2)^2)^{\frac 12} \: (\int\limits_D |\nabla v|^2 dx)^{\frac 12}\qquad \forall v\in H^1_0(D)$$
Thus $J_K$ achieves its global minimum, which is the unique $H^1_0(D)$ solution of (\ref{curvature zero}). The same arguments, applied to $K_\Phi= K e^{2\Phi}$, where $\Phi$ is the harmonic extension of some boundary data $\varphi$, lead to the following
\begin{teor} \label{positive}
Let $K\leq 0$ and $\int\limits_D K^2 (1-|x|^2)^2 dx < +\infty$. \\Given a smooth boundary data $\varphi$, (\ref{curvature zero})has a unique  solution which takes the boundary data $\varphi$ and which writes as $u=v+\Phi, v\in H^1_0(D)$. \\In particular, $K$ is the curvature of the conformal metric $g= e^{2(v+\Phi)}g_e$.
\end{teor} 
To get a \textsl{complete} conformal metric with curvature $K$, one can build, following \cite{LN}, a sequence $u_n$ of solutions of (\ref{curvature zero}) taking $\varphi\equiv n$ and try to show that it converges to a solution $u$ of (\ref{curvature zero})
such that $u(x)\rightarrow +\infty$ suitably fast as $|x|\rightarrow 1$. We don't pursue the details.\\
A more natural approach to find a complete conformal metric with curvature $K$, is  to look for a \emph{bounded} $ C^2$ solution of the equation
\begin{equation} \label{curvature -1}
\Delta_{\mathcal H} u  + 1 + K e^{2u} = 0
\end{equation}
where  $\Delta_{\mathcal H} $ denotes the hyperbolic laplacian (notice that  
 solutions $u$ of (\ref{curvature -1}) and $v$ of (\ref{curvature zero}) are simply related: $v-u= \log \frac{2}{1-|x|^2}$).  We recall the following pioneering result (\cite{AMo}, see also \cite{BK}) 
\begin{teor}(Aviles-McOwen) \label{AO} 
Let $K\in C^\infty(D)$, $K\leq 0$ in $D$ and such that $-\frac{1}{c}\leq K \leq -c $ in $\{c \leq |x| < 1\}$ for some $c\in (0,1)$. Then there is a unique metric conformal and uniformly equivalent to the hyperbolic metric having $K$ as its Gaussian curvature.
\end{teor}
We end this section with a result which might provide complete conformal metrics with prescribed nonpositive gaussian curvature. Given a conformal metric $g$ on the disc, let us denote by $K_g$ its curvature. Given $K$, $e^{2u}g$ is a conformal metric with curvature $K$ if 
   $u\in C^2(D)$  satisfies the equation
\begin{equation} \label{curvature}
\Delta_g u -K_g + K e^{2u} = 0
\end{equation}
If, in addition, $u$ is bounded, then $e^{2u}g$ is quasi isometric to $g$. In this case, if $g$ is complete then   $e^{2u}g$ is complete as well.
\begin{teor} \label{AOms} 
Let $g\leq c g_h$ be a conformal metric. Let $K=K_g+H$ be  nonpositive in $D$. \\
Assume  $H \in L^2(D, \mu_g)$ . Then (\ref{curvature}) has a solution in  $H^1_0$.
\end{teor}
\begin{proof} Solutions for (\ref{curvature}) can be  obtained as critical points of the functional
$$J_K(v)= \frac 12 \int\limits_D |\nabla v|^2 dx - \int\limits_D H v dV_g - \frac 12 \int\limits_D K (e^{2v}-2v-1) dV_g \qquad v\in H^1_0(D)$$
The assumption on $g$ implies $\lambda_1(g):=S_2(D,g)>0$ and hence 
$$J_K(v) \geq \frac 12 \int\limits_D |\nabla v|^2 dx - \frac{1}{\sqrt{\lambda_1(g)}}(\int\limits_D H^2 dV_g)^{\frac 12}(\int\limits_D |\nabla v|^2 dx)^{\frac 12}\qquad \forall v\in H^1_0(D)$$
Thus $J_K$ is  a  (possibly infinite somewhere)  convex coercive functional in $H^1_0(D)$. 
By Fatou's Lemma it is also weakly lower semicontinuous, and hence it achieves its infimum at some $\underline v$.\\
Notice that $J_K(\underline{v}+t\varphi) <+\infty$ for all $\varphi \in C_0^\infty(D)$ because
$$\int\limits_{supp (\varphi)} (-K) (e^{2(\underline v+t\varphi)}-2(\underline v+t\varphi)-1)dV_g \leq
\sup_{supp (\varphi)} (-K) \int\limits_D  (e^{2\underline v}-2 \underline v-1) dV_g < +\infty$$
by Trudinger exponential integrability.
Hence
$$0=\frac{d}{dt}J_K(\underline v+t\varphi)_{|_{t=0}}=\int\limits_D \nabla_g \underline v \nabla_g \varphi -
(K-K_g)\varphi - K (e^{2 \underline v}\varphi-\varphi) dV_g$$
i.e. $\underline v$ solves (\ref{curvature}).
\end{proof}
\begin{rmk} In particular, following \cite{AMo}, one can take  $K=f+H$, $f\in L^2(D, \mu_h)$ and $H\leq 0$ bounded and bounded away from zero around $\partial D$. 
\end{rmk}
\begin{rmk} The above result slightly improves a result by D.M. Duc \cite{D}, where, in addition, conditions are given to insure the metric is complete.
\end{rmk}
\setcounter{equation}{0}

\section{Appendix} \label{A}
We present an example of a domain for which $\omega(\Omega) < +\infty$ and $\lambda_1(\Omega)= 0$. Let
$$\Omega = \R^2 \setminus \bigcup_{n, m \in \Z}D_{n,m} \qquad D_{n,m}= D_{r_{n,m}}(n,m)\qquad \log \frac{1}{r_{n,m}} = 2^{|n|+|m|}$${}{}{}{}
We are going to exibit a sequence $u_k\in H^1_0(\Omega)$ such that 
$$ \sup_k \int\limits_\Omega |\nabla u_k|^2 < \infty \qquad \int\limits_\Omega u_k^2 \rightarrow_k +\infty$$
Let $\psi_k \in C_0^\infty (D_{3k}, [0,1]), \quad \psi_k \equiv 1$ in $D_{ k}$, be  radial  with $|\nabla\psi_k| \leq \frac 1k$, so that 
$$\int\limits_{\R^2} | \nabla \psi_k|^2 \leq 8 \pi \quad \emph{and} \quad \int\limits_{\R^2}|\psi_k|^2 \geq \pi k^2$$
Let $\varphi_{\epsilon}(x)= 2(1-\frac{\log |x|}{\log \epsilon})$ in $A_\epsilon := \{\epsilon \leq |x| \leq \sqrt \epsilon \}$ and
$\varphi_\epsilon \equiv 0$ in $|x|\leq \epsilon$, so that
$$\int\limits_{|x|\leq \sqrt\epsilon} |\nabla \varphi_\epsilon|^2 \leq -\frac{4\pi}{\log \epsilon}  \quad \emph{and} \quad \int\limits_{|x| \leq \sqrt \epsilon}|\varphi_\epsilon|^2 \leq \epsilon \pi$$
Finally, let $\varphi = \varphi_{r_{n,m}}(x- (n,m))$ in $D_{\sqrt{r_{n,m}}}(n,m)$, $\varphi = 1$ elsewhere, and let
$$u_k (x)= \min \{\varphi(x), \psi_k(x)\}$$
so that $u_k \in H^1_0(\Omega)$ and 
$$\int\limits_{\Omega}|\nabla u_k|^2 \leq 8  \pi + 4 \pi \sum_{n,m}\frac{1}{2^{|n|+|m|}} \leq 44 \pi \quad \emph{and} \quad \int\limits_{\Omega}| u_k|^2 \geq \pi k^2 - \pi \sum_{n,m} r_{n,m}$$

\end{document}